\documentclass[11pt,a4paper]{article}
\usepackage[left=27mm,right=27mm,top=23mm,bottom=25mm]{geometry}
\usepackage[T1]{fontenc}
\usepackage{lmodern,amsmath,amssymb,amsthm,graphicx,booktabs,xcolor}
\usepackage[colorlinks=true,linkcolor=teal,citecolor=teal,urlcolor=teal]{hyperref}
\numberwithin{equation}{section}
\newtheorem{theorem}{Theorem}[section]
\newtheorem{lemma}[theorem]{Lemma}
\newtheorem{proposition}[theorem]{Proposition}

\newcommand{\area}{\operatorname{area}}
\newcommand{\conv}{\operatorname{conv}}
\title{The smallest vertex triangle in a fixed polygon}
\author{Dawid Trela\\[2pt]\small Doctoral School, War Studies University, Warsaw, Poland\\
\small \href{mailto:dawidmtrela@gmail.com}{dawidmtrela@gmail.com}\quad
\href{https://orcid.org/0000-0001-9781-6425}{ORCID 0000-0001-9781-6425}}
\date{}
\begin{document}
\maketitle
\begin{abstract}
Among the triangles determined by $n$ points in convex position in the unit square, one has area at most $32/n^3$. Equality is possible exactly when $4$ divides $n$, and then the maximizing point set is unique. We prove a fixed-container theorem for a class of convex $s$-gons that includes affinely regular polygons and a one-parameter family of pairwise affinely inequivalent hexagons without a centre of symmetry. If the container has equal ears of area $a$ and the tips of its directed side vectors are in convex position, the sharp bound at $n=sk$ is $a/k^3$. Its equality conditions determine every side vector of the maximizing polygon. The proof uses the area enclosed by the normalized side directions. We compare the result with the area-normalized smallest-ear theorem of Gronchi and Longinetti and give consequences for unrestricted Heilbronn configurations with 13 and 15 points.
\end{abstract}

\section{An octagon in a square}
\begin{figure}[!htb]
\centering\includegraphics[width=.93\linewidth]{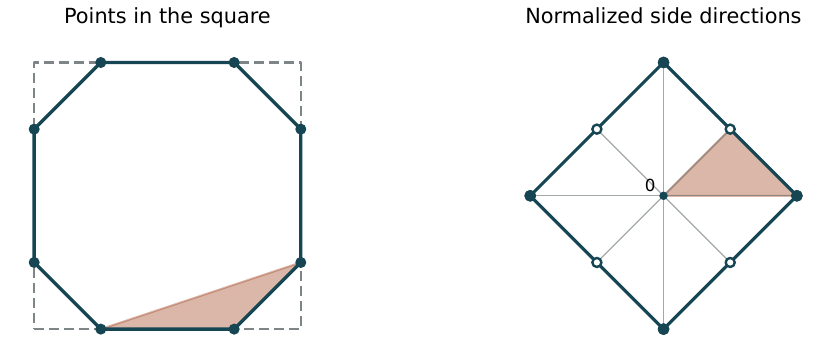}
\caption{The optimal octagon in the unit square and the tips of its eight normalized edge vectors. The shaded triangle on the right has area equal to half the determinant of two successive directions.}\label{fig:square}
\end{figure}

Mark the two quarter points of each side of the unit square. Take the two marks on the bottom side and the nearest mark on the right. Their triangle has base $1/2$, height $1/4$ and area $1/16$ (Figure~\ref{fig:square}). No triangle on the eight marks is smaller. Could eight other points in convex position do better? In fact even a tie forces the eight marks back to these positions.

For a finite set $P$ of at least three points, let $m(P)$ be the least area of a triangle determined by three of its points. For a compact convex region $K$ with interior, let
\[
 C_n(K)=\sup\{m(P):P\subset K,\ |P|=n,\ \text{every point of }P\text{ is a vertex of }\conv P\}.
\]
We use $|P|$ for the cardinality of a finite point set and $\area(K)$ for area. The supremum is attained. A small regular $n$-gon in $K$ gives a positive floor. Label the sets in a maximizing sequence in counterclockwise order and pass to a convergent subsequence in $K^n$. Every ordered triple in the limit still has strictly positive signed area. Thus each boundary edge has all the other vertices strictly on its left, and the limiting set is admissible.

\begin{theorem}\label{th:square}
In the unit square $Q$, $C_n(Q)\leq 32/n^3$ for every $n\geq3$. Equality holds exactly at multiples of four. For $n=4k$, the maximizing set is unique and its least triangle has area $1/(2k^3)$.
\end{theorem}

This problem has two close relatives. Barequet~\cite{Barequet} obtained cubic decay for points on a convex decreasing chain in a square. Closing the chain brings in its entire turn and makes it possible to identify the sharp constant. Ismailescu, Kim and Wang~\cite{IKW} allow the shape of the outer $n$-gon to vary and measure the area lost when a vertex is removed as a fraction of the polygon's area. Their fraction is governed by the smallest ear, the triangle on three consecutive vertices. Gronchi and Longinetti~\cite{GL} identified the maximizers of this ratio. Section~\ref{sec:affine} recalls the sharp value and explains its connection with Fejes Tóth's affine-perimeter inequality. Here the outer figure is prescribed, and two containers of the same area can give different answers.

Convex polygons with equal ears have an established discrete affine geometry~\cite{Craizer}. Qi et al.~\cite{Qi} study area minimization under equal consecutive-triangle constraints. For regular containers the midpoint parabolas of affine-perimeter theory already suggest a construction~\cite{Fejes,Morin}. For a finite number of vertices, the question is which side vectors an extremizer must have. Our proof gives them explicitly and works for a container without affine regularity. Vinciguerra's thesis seminar~\cite{Vinciguerra} also treats local optima and decreasing configurations. Section~\ref{sec:free} gives an application to the unrestricted Heilbronn problem.

\section{A perimeter and an area of directions}
Throughout, a convex polygon is listed counterclockwise and each listed vertex is extreme. For the vertices $p_i$ of a candidate polygon put $f_i=p_{i+1}-p_i$, with cyclic subscripts. The ear at $p_i$ has twice its area equal to $\det(f_{i-1},f_i)>0$.

\begin{lemma}\label{lem:ear}
At least one of the smallest vertex triangles of a convex polygon is an ear.
\end{lemma}
\begin{proof}
Fix a chord $AB$ and one of its boundary arcs with intermediate vertices. On that arc the distances to the line $AB$ increase and then decrease, allowing ties. Indeed, a line parallel to the chord meets the convex boundary in at most two points, except at a supporting segment. At each positive height the part of the arc above that height is consequently an interval. The least height among the intermediate vertices is attained at the first or the last.

Start with a smallest triangle $ABC$. Keeping $AB$ fixed, replace $C$ by a vertex of least height on its arc. This gives another minimum triangle with a polygon edge as one side. Use that edge as the new base. On its complementary arc the least height is again attained next to an endpoint. The resulting triangle is an ear.
\end{proof}

Choose a compact convex set $D$ containing the origin in its interior. Measure a vector by the factor by which one must dilate $D$ for its boundary to reach that vector. More precisely, its gauge is $\rho_D(x)=\inf\{t>0:x\in tD\}$. For a counterclockwise polygon $P$ write
\[
 \ell_i=\rho_D(f_i),\qquad u_i=f_i/\ell_i,\qquad
 p_D(P)=\sum_i\ell_i,\qquad U=\conv\{u_0,\ldots,u_{n-1}\}.
\]
In the square example $D$ is the diamond with vertices $(1,0),(0,1),(-1,0),(0,-1)$. Here $\ell_i$ is the sum of the absolute coordinate increments. Figure~\ref{fig:square} distinguishes the polygon in the container from the polygon of edge directions in $D$. The construction also works for an asymmetric $D$. In that case reversing a vector can change its length, so the orientation of the boundary is part of the definition.

\begin{theorem}\label{th:gauge}
If an $n$-vertex convex polygon $P$ is contained in a convex polygon $K$, then
\begin{equation}\label{eq:bound}
 n^3 m(P)\ \leq\ p_D(P)^2\,\area(U)
       \ \leq\ p_D(K)^2\,\area(D).
\end{equation}
Equality throughout requires all the lengths $\ell_i$ to agree, all determinants $\det(u_{i-1},u_i)$ to agree, $U=D$, and $p_D(P)=p_D(K)$.
\end{theorem}
\begin{proof}
The directions of the $f_i$ turn once counterclockwise. The $u_i$ have that same order on $\partial D$. Their positive linear combination $\sum_i\ell_i u_i=0$ and their span is the plane. A supporting line through the origin would then contain them all, so the origin is interior to $U$. Each $u_i$ lies on the boundary of $U$, since it lies on a support line of $D$. In their radial order these points therefore trace $\partial U$, including any points along a straight side. Joining each side of this circuit to the origin gives
\begin{equation}\label{eq:fan}
 q_i=\det(u_{i-1},u_i)>0,\qquad \sum_i q_i=2\area(U).
\end{equation}

For every $i$, $2m(P)\leq\ell_{i-1}\ell_i q_i$. Multiply these $n$ inequalities. Each length occurs twice. Applying the arithmetic--geometric mean inequality to the $n$ lengths and then to the $n$ determinants gives
\[
 2m(P)\leq
 \left(\prod_i\ell_i\right)^{2/n}\left(\prod_i q_i\right)^{1/n}
 \leq\left(\frac{p_D(P)}{n}\right)^2\frac{2\area(U)}{n}.
\]
The two equality conditions in these means give the stated conditions on the lengths and determinants.

Finally $p_D(P)\leq p_D(K)$. Intersect $K$ successively with the closed half-planes supporting the sides of $P$. A cut that removes a cap replaces a directed chain by the chord having the same vector sum. Convexity of $D$ gives $\rho_D(x+y)\leq\rho_D(x)+\rho_D(y)$, by taking the length-weighted convex combination of the normalized vectors. Hence a cut cannot increase the perimeter. After the last cut the polygon is $P$. Since $U\subset D$, its area is at most that of $D$, with equality only for $U=D$. Indeed a proper closed convex subset misses a cap of positive area.
\end{proof}

\section{Which polygons attain the bound?}\label{sec:sharp}
Let $K$ have $s$ vertices $v_0,\ldots,v_{s-1}$, and write $e_r=v_{r+1}-v_r$. Translate the tail of each directed edge to the origin. The next theorem uses the hull $D$ of the resulting $s$ tips $e_r$. Its hypotheses are finite determinant and convexity tests on the given container.

\begin{theorem}\label{th:exact}
Assume that the $e_r$ are the distinct successive vertices of $D$ in counterclockwise order and that $\det(e_{r-1},e_r)=2a$ is independent of $r$. Then
\begin{equation}\label{eq:exact}
 C_n(K)\leq \frac{s^3a}{n^3},\qquad
 C_{sk}(K)=\frac{a}{k^3}\quad(k=1,2,\ldots).
\end{equation}
Equality in the upper bound holds precisely when $s\mid n$. At each such $n$ there is exactly one maximizing point set in the fixed container.
\end{theorem}
\begin{proof}
The origin is in the interior of $D$ because the $e_r$ span the plane and sum to zero. As every $e_r$ is on $\partial D$, each side of $K$ has gauge length one. Joining the origin to the sides of $D$ gives $s$ triangles, each of area $a$. Thus $p_D(K)=s$ and $\area(D)=sa$. The upper bound follows from Theorem~\ref{th:gauge}.

\smallskip\noindent\emph{What equality forces.}
Suppose the upper bound is attained. All the normalized directions together have hull $D$, so every corner $e_r$ must occur. The directions partition each side of $D$ in their cyclic order. On the side from $u=e_{r-1}$ to $v=e_r$ set $w(t)=(1-t)u+tv$. Bilinearity gives
\begin{equation}\label{eq:line}
 \det(w(t),w(t'))=(t'-t)\det(u,v)=2a(t'-t).
\end{equation}
By~\eqref{eq:fan}, the common $q_i$ is $2sa/n$. Every step on each side therefore has parameter increment $s/n$. Traversing a side of parameter length one takes exactly $n/s$ steps. Write $k=n/s$ and $h=1/k$. The common gauge length is $s/n=h$. Up to the initial index, the entire list of side vectors is now forced. It runs through $h\,w(jh)$ for $j=0,\ldots,k-1$ on each successive side of $D$.

\smallskip\noindent\emph{Where to put those edges.}
Near the corner $v_r$ put $k$ points, one for each $t=jh$ with $0\leq j<k$:
\begin{equation}\label{eq:corner}
 z_{r,j}=v_r-\frac{(1-t)(1-t-h)}{2}\,e_{r-1}
                  +\frac{t(t+h)}{2}\,e_r .
\end{equation}
Denote the two nonnegative coefficients by $A,B$. They obey $A\leq(1-t)/2$ and $B\leq t/2$. The point is the convex combination of $v_{r-1},v_r,v_{r+1}$ with weights $A,1-A-B,B$, and lies in $K$. Subtracting consecutive points in one block gives $h[(1-t-h)e_{r-1}+(t+h)e_r]$. Across a junction the points are $v_r+(1-h)e_r/2$ and $v_r+(1+h)e_r/2$. Their difference is $he_r$. Thus the closed circuit has exactly the prescribed edge vectors, also when $k=1$.

Their directions turn strictly through one revolution, every step being less than $\pi$. To see that this circuit is convex, choose any of its directed edges as a base. As the remaining edges are followed, their determinants with the base are first nonnegative and then nonpositive. The first is positive, the last negative, and the height returns to zero on closing. All other vertices are strictly to the left of the base line. Hence all $sk$ points are distinct extreme vertices. By~\eqref{eq:line}, two successive edges have determinant $h^2(2ah)=2ah^3$. Every ear has area $a/k^3$, which is also the smallest vertex-triangle area by Lemma~\ref{lem:ear}.

\smallskip\noindent\emph{Why the position is unique.}
Polygons with the same cyclic list of directed side vectors differ by a translation. Our construction contains a segment centred on each side of $K$, of relative length $h$. If its translate by $b$ still lies in $K$, every outward side normal $\nu_r$ satisfies $\nu_r\cdot b\leq0$. A nonzero such $b$ would keep the entire ray $x+tb$ in $K$ for $x\in K$ and $t\geq0$. Boundedness forces $b=0$.
\end{proof}

For $k=2$ formula~\eqref{eq:corner} places the new points a quarter of a side from each vertex, exactly one step of Chaikin's corner-cutting construction~\cite{Chaikin}. For general $k$ the same points are the intersections of tangents at successive subdivision nodes on the quadratic arcs used in Section~\ref{sec:affine}. Indeed the two tangent lines contain the same point $\gamma(t)+h\gamma'(t)/2=\gamma(t+h)-h\gamma'(t+h)/2$, by the quadratic Taylor formula. The minimum-area condition thus characterizes a familiar construction.

A regular $s$-gon satisfies both hypotheses by rotational symmetry. The conditions survive an invertible affine map, with areas multiplied by the absolute determinant of its linear part (reverse the boundary order if needed). In particular for a triangle $T$ and a parallelogram $S$ we obtain
\[
 C_{3k}(T)=\area(T)/k^3,\qquad C_{4k}(S)=\area(S)/(2k^3).
\]
The latter proves Theorem~\ref{th:square}. For example, in a triangle of area one the twelve-point value is $1/64$, while a parallelogram of the same area gives $1/54$. For any $K$ of Theorem~\ref{th:exact}, delete vertices of its optimizer at the next multiple of $s$. Deletion cannot reduce the minimum, so
\[
 \frac{a}{\lceil n/s\rceil^3}\leq C_n(K)\leq\frac{s^3a}{n^3},\qquad
 \lim_{n\to\infty}n^3 C_n(K)=s^3a .
\]

\subsection*{A family of hexagons without a centre of symmetry}
Fix $0<u<v<2u$. Starting from the origin, take the directed steps
\[
 (u,0),\quad (v,v),\quad (0,u),\quad (-v,0),\quad (-u,-u),\quad (0,-v).
\]
They sum to zero and their angles make one turn, giving a convex hexagon $K$. Consecutive steps have determinant $uv$, so $a=uv/2$. To check the other hypothesis, join their tips in order. The circuit is simple since each chord stays in its own angular sector. Its successive turning determinants alternate between $v(2u-v)$ and $u(2v-u)$, both positive. Thus every tip is an extreme point. Theorem~\ref{th:exact} gives $C_{6k}(K)=uv/(2k^3)$.

\begin{figure}[!htb]
\centering\includegraphics[width=.95\linewidth]{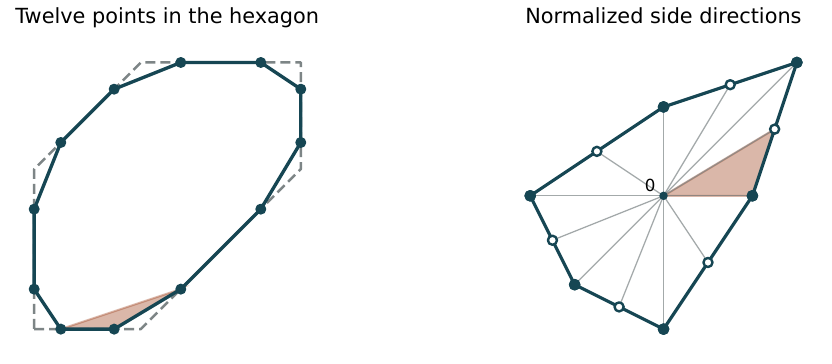}
\caption{Left, the twelve points obtained by subdividing the six corners of the given hexagon. Right, the six directed side vectors with their tails at the origin. The half steps on their convex hull determine every edge of the optimizer.}\label{fig:hexagon}
\end{figure}

Each of the three parallel pairs has longer-to-shorter side ratio $v/u$. This ratio is preserved by invertible affine maps. Fixing $u=1$ and varying $v$ between one and two therefore gives pairwise affinely inequivalent hexagons. Their unequal opposite sides preclude central symmetry. For $(u,v)=(2,3)$ the vertices are $(0,0),(2,0),(5,3),(5,5),(2,5),(0,3)$. Its unique twelve-point optimizer has minimum area $3/8$ (Figure~\ref{fig:hexagon}), and its eighteen-point optimizer has minimum area $1/9$.

\section{The smallest ear relative to the area of its hull}\label{sec:affine}
Gronchi and Longinetti~\cite{GL} proved that the ratio of the smallest ear to the polygon area is maximized exactly by affinely regular polygons. The next proposition is their area-normalized result, expressed also in the notation of Ismailescu, Kim and Wang~\cite[Section~4]{IKW}. We give the short affine-perimeter argument because it shows precisely which part of the regular-container bound is already available from classical theory.

\begin{proposition}\label{prop:normalized}
For a convex $n$-gon $P$, $n\geq3$,
\[
 \frac{m(P)}{\area(P)}\leq\frac4n\sin^2(\pi/n).
\]
Equality holds precisely for affine images of a regular $n$-gon. In the notation of~\cite{IKW}, this gives $1-f_n(n-1)=(4/n)\sin^2(\pi/n)$ for $n\geq4$.
\end{proposition}
\begin{proof}
Let $A_i$ be the ear areas. At a vertex $p_i$ join its two adjacent side midpoints $M_-,M_+$ by the quadratic arc $(1-t)^2M_-+2t(1-t)p_i+t^2M_+$, $0\leq t\leq1$. It lies in their corner triangle. Its velocity is $(1-t)f_{i-1}+tf_i$, so the determinant of velocity and acceleration is $2A_i$. The arcs join with matching tangents, and these tangents turn strictly through one revolution. Relative to any fixed tangent the signed height first increases, then decreases to zero. Thus the arcs bound a convex body in $P$. In the normalization $\int |\det(\gamma',\gamma'')|^{1/3}\,dt$, its affine perimeter is $\lambda=\sum_i(2A_i)^{1/3}$.

Fejes Tóth's inequality for a convex body in an $n$-gon~\cite[Ch.~II, \S6, (1)]{Fejes} now gives
\[
 2n^3m(P)\leq\lambda^3\leq8\,\area(P)\,n^2\sin^2(\pi/n).
\]
Equality in the right inequality requires the enclosing $P$ to be affinely regular. In a regular $n$-gon of circumradius $R$, an ear has base $2R\sin(2\pi/n)$ and height $R(1-\cos(2\pi/n))$. Dividing its area by $nR^2\sin(2\pi/n)/2$ gives the displayed ratio. It is preserved by affine maps. Finally, a largest inscribed $(n-1)$-gon can be chosen by deleting one vertex~\cite{IKW}. The area lost is the smallest ear area, equal to $m(P)$ by Lemma~\ref{lem:ear}.
\end{proof}

In a regular container the maximum affine perimeter is realized by the midpoint parabolas just used~\cite[Appendix A, Lemmas 16--17]{Morin}. Applying the same rounding to each candidate already gives the upper bound of Theorem~\ref{th:exact} in this special case. Its elementary proof also supplies the allowed cardinalities, the individual edge vectors and the unique placement for the larger class of fixed containers.

\section{Convex subsets of unrestricted configurations}\label{sec:free}
Let $H_n$ be the unrestricted Heilbronn optimum in the unit square. A configuration with $h$ points in convex position satisfies $h^3m(P)\leq32$ by Theorem~\ref{th:square}, since every triangle on the subset is a triangle of $P$. For 13 and 15 points, respectively, the configurations of Karpov and Sudermann-Merx in~\cite[Appendix B]{Sudermann}, rounded to a grid of mesh $1/1000$, give
\[
 H_{13}\geq\frac{53341}{2\,000\,000}>\frac{32}{11^3},\qquad
 H_{15}\geq\frac{41853}{2\,000\,000}>\frac{32}{12^3}.
\]
The supplement lists the coordinates and gives exact integer checks for each triple. Thus every unrestricted maximizer has at most ten points in convex position for $n=13$, and at most eleven for $n=15$. At least three and four points, respectively, lie strictly inside their hulls. The exact unrestricted values of $H_{13}$ and $H_{15}$ are not determined by these exclusions. Another question left by the square theorem is the value of $C_n(Q)$ between successive multiples of four.

\section*{Supplementary material}
The project in Lean 4 proves the square upper bound and its attainment for every $4k$ in a model of a counterclockwise polygonal boundary. Each edge is given by its coordinates, a positive length and an ordered position on the direction diamond. A convex polygon has this presentation after cyclic relabelling: assign an axial direction to the quarter starting there, then read the four consecutive blocks in order, allowing empty blocks. Conversely a positive signed-area floor on all ordered triples places every other vertex strictly to the left of each cyclic boundary edge. The declaration \texttt{Heilbronn.square\_optimum} quantifies over every competitor in this representation.

The project also checks the finite product estimate in Theorem~\ref{th:gauge}, the subdivision identities, the determinant and orientation hypotheses across the hexagon family, the side count forced by equal increments, and both unrestricted grid witnesses. The geometric gauge argument, the general equality classification and Proposition~\ref{prop:normalized} are proved in the text. The Lean proof and the plotting and exact-arithmetic scripts were developed with OpenAI ChatGPT (Codex). The drawings are deterministic plots of formula~\eqref{eq:corner}. A separate note records a matching of interior points to distinct minimum triangles and two structural results from the author's earlier work~\cite{Trela1,Trela2}.

\section*{Declaration of generative AI and AI-assisted technologies}
OpenAI ChatGPT (Codex) assisted with the development and checking of proof arguments, the literature search, drafting and editing, the Lean formalization, and the scripts for the numerical checks and figures. Figures~\ref{fig:square} and~\ref{fig:hexagon} plot the stated coordinates using the supplied script.

\end{document}